\documentclass{amsart}

\usepackage{amssymb}
\usepackage[all]{xy}
\usepackage{hyperref}

\usepackage{enumitem}   

\setlist[enumerate]{itemsep=.2em,topsep=.2em,leftmargin=1.25em,itemindent=2.0em}

\newtheorem{thm}{Theorem}

\newtheorem{lem}[thm]{Lemma}
\newtheorem{cor}[thm]{Corollary}

\newtheorem{prop}[thm]{Proposition}

\newtheorem{conj}[thm]{Conjecture}
   
\theoremstyle{definition}
\newtheorem{defn}[thm]{Definition}

\newtheorem{say}[thm]{}

\newtheorem{rems}[thm]{Remarks} 

\newtheorem*{ack}{Acknowledgments}      

\newtheorem{defn-thm}[thm]{Definition--Theorem}  
\newtheorem{defn-lem}[thm]{Definition--Lemma}  

\newtheorem{main-exmp}[thm]{Main Example}
\newtheorem{baby-exmp}[thm]{Baby Example}
\newtheorem*{obs-nn}{Observation}
\newtheorem*{assertion-nn}{Assertion}   

\theoremstyle{remark}

\renewcommand{\o}[0]{{\mathcal O}} 

\newcommand{\p}[0]{{\mathbb P}}

\newcommand{\map}[0]{\dasharrow}
\newcommand{\qtq}[1]{\quad\mbox{#1}\quad}

\newcommand{\br}[0]{\operatorname{Br}}

\newcommand{\grass}[0]{\operatorname{Grass}}

\newcommand{\univ}[0]{\operatorname{Univ}}

\newcommand{\bire}[0]{\operatorname{\stackrel{bir}{\sim}}}

\newcommand{\ind}[0]{\operatorname{ind}}

\def\into{\DOTSB\lhook\joinrel\to}

\def\loccoh#1.#2.#3.#4.{H^{#1}_{#2}(#3,#4)}

\DeclareMathAlphabet{\mathchanc}{OT1}{pzc}%
                                {m}{it}

\usepackage[all]{xy}\xyoption{dvips}

\begin{document}
\bibliographystyle{amsalpha}


 \title[Severi--Brauer varieties]{Birational equivalence of Severi--Brauer varieties}
 \author{J\'anos Koll\'ar}

 \begin{abstract} We prove Amitsur's conjecture for Severi--Brauer varieties whose index is not a prime power.
       \end{abstract}

  \maketitle

   Let $P$ be  a Severi--Brauer variety over a field $k$.
 Amitsur  \cite{MR70624} proved that if $P$ is birational to another  Severi--Brauer variety $Q$, then they 
 generate the same subgroup of $\br(k)$, and conjectured that the converse also holds, provided $P, Q$  have the same dimension.
The conjecture was popularized by Serre in 
 \cite[Sec.X.6]{MR150130, MR554237}.

 Roquette  \cite{MR166215} developed an inductive  method, and settled  the conjecture when $P$ splits over a solvable Galois extension of degree $\dim P+1$.
 Krashen   \cite {MR2029172, MR2363485} used the method to show, among others,  that if the conjecture holds for both $P_1, P_2$, and $\dim P_1+1, \dim P_2+1$ are relatively prime,
 then it also holds for their (tensor or Brauer)  product $P_1\otimes P_2$; see Definition~\ref{tensor.pr.defn}.
 Combined with \cite{tregub}, this implies that the conjecture holds if 2 is a primitive root modulo every prime factor of  $\dim P+1$.

 \smallskip
 We show that the conjecture holds for  $P_1\otimes P_2$, even if we do not know what happens for  $P_1$ and $P_2$ individually.
 This leads to the following.

 \begin{thm} \label{stab.bir.thm} 
 Let $P,  Q$ be  Severi--Brauer varieties  of the same dimension over a  field $k$. Assume that $\dim P+1$ has at least 2 distinct prime factors.
Then $P, Q$ generate the same subgroup of $\br(k)$ iff  $P$ is birationally equivalent  to $Q$.
 \end{thm}

 This turns out to be a  straightforward consequence of the following birational description of the tensor product of Severi--Brauer varieties.
 
 \begin{thm}\label{bir.of.prod.thm}  Let $P_1,P_2$ be  Severi--Brauer varieties over a field  $k$, and assume that
 $(\dim P_1+1, \dim P_2+1)=1$. Then there is a  birational equivalence
   $$
   P_1\otimes P_2\bire P_1\times P_2\times \p^r,\qtq{where}
   r=\dim P_1\cdot \dim P_2.
   \eqno{(\ref{bir.of.prod.thm}.1)}
   $$
 \end{thm}

 \begin{rems} \label{bir.of.prod.thm.rems}{\ }

   (\ref{bir.of.prod.thm.rems}.1) The method of \cite{MR70624} produces 
   birational equivalences of the form
   $$
   (P_1\otimes P_2)\times \p^{r_1}\bire P_1\times P_2\times \p^{r_2}
   $$
   for some $r_1, r_2>0$, see (\ref{conj.conj}.4).
   The key assertion of  Theorem~\ref{bir.of.prod.thm} is that we can achieve $r_1=0$.
   
   (\ref{bir.of.prod.thm.rems}.2) The birational equivalence (\ref{bir.of.prod.thm}.1) usually fails without the realtively prime assumption.
   For example, if $P_1=P_2=P$  and $\dim P+1$ is even,  then
   $P\otimes P$ and $P\times P$ are not  even stably birational to each other.

  (\ref{bir.of.prod.thm.rems}.3) There is no general birational product formula for 
   $P_1\otimes P_2$ that involves only the $P_i$ and their dimensions.
   For example,  although $P$ is birationally equivalent to its dual $P^\vee$, the tensor product  $P\otimes P^\vee$ is birationally trivial, but
   $P\otimes P$ is usually not. See, however, Lemma~\ref{prod.to.power.lem}.
  \end{rems}

 Note that in Theorem~\ref{bir.of.prod.thm},
 the subgroup of $\br(k)$ generated by
 $P_1\otimes P_2$ is the same as
 the subgroup  generated by
 $P_1$ and $P_2$. This suggests the following generalization.

 \begin{conj} \label{conj.conj} Let $X=\prod_{i\in I} P_i$ and  $Y=\prod_{j\in J} Q_j$ be (ordinary) products of  Severi--Brauer varieties.
   The following are equivalent.
   \begin{enumerate}
   \item The subgroups  $\langle P_i: i\in I\rangle\subset \br(k)$ and
     $\langle Q_j: j\in J\rangle\subset \br(k)$ agree.
   \item  $X\times \p^{\dim Y}\bire   Y\times \p^{\dim X}$.
       \end{enumerate}
\noindent  If  $\dim X=\dim Y$, then these are also  equivalent to
\begin{enumerate}\setcounter{enumi}{2}
\item $X \bire   Y$.
   \end{enumerate}
 \end{conj}

 {\it Known cases.}  Using Corollary~\ref{BR.BC.COR}  we see that
 the pull-back of each $Q_j$  to $X$ is trivial. Similarly,
 the pull-back of each $P_i$  to $Y$ is trivial. This gives
 birational equivalences
 $$
 X\times \p^{\dim Y}\bire  X\times Y \bire  Y\times \p^{\dim X}.
 \eqno{(\ref{conj.conj}.4)}
 $$
 Thus (\ref{conj.conj}.1) and (\ref{conj.conj}.2) are equivalent.

 (\ref{conj.conj}.3)  $\Rightarrow$ (\ref{conj.conj}.2) is clear, and 
repeated application of Theorem~\ref{bir.of.prod.thm} gives many cases of the converse.
 
 Some cases where $P_i, Q_j$ have dimension $\leq 2$ are treated in 
\cite{k-conic, MR2439423}.

 \subsection*{Background on Severi--Brauer varieties}{\ }

 We collect  basic results about Severi--Brauer varieties that we need. The algebraic versions are discussed in several books \cite{MR595, MR227205, MR2266528}.
 A geometric treatment is given in \cite{k-sb}.

 \begin{defn} \label{sb.defn} Let $k$ be a field with algebraic closure $\bar k$. A {\it Severi--Brauer variety}  is a $k$-variety $P$  such that $P_{\bar k}\cong \p^{\dim X}_{\bar k}$.

   $Q\subset P$ is a  {\it twisted linear} subvariety if
   $Q_{\bar k}\subset P_{\bar k}\cong \p^{\dim X}_{\bar k}$ is a linear subvariety.

   The {\it index} of $P$ is the gcd of the degrees of 0-cycles  on $P$.
It is also the minimum degree of a (nonzero) effective 0-cycle  on $P$.
   The index divides $\dim P+1$; it is denoted by $\ind(P_k)$.
   \end{defn}
 
 \begin{lem} \label{P.min.lem}
   Let $P$ be  a Severi--Brauer variety.
   All minimal twisted linear subvarieties are isomorphic to each other; denote these by   $P^{\rm min}$.

 Then   $\dim P^{\rm min}=\ind (P)$,  and $P$ is birationally equivalent to
   $P^{\rm min}\times \p^r$ where $r=\dim P-\dim P^{\rm min}$.
   \qed
 \end{lem}

\begin{defn}[Tensor or Brauer product]  \label{tensor.pr.defn}
   Let $P,  Q$ be  Severi--Brauer varieties   over a  field $k$.
   There is a unique Severi--Brauer variety---denoted by
   $P\otimes Q$---and an embedding
   $P\times Q\into P\otimes Q$, that is isomorphic to the Segre embedding over $\bar k$.
\end{defn}

\begin{defn}[Brauer group] \label{br.gr.defn}
  Two Severi--Brauer varieties $P, Q$ are   {\it Brauer equivalent} iff $P^{\rm min}\cong Q^{\rm min}$.

  The tensor product defines a group structure on the set of
  Brauer equivalence classes over $k$. It is called the
  {\it Brauer group} of $k$, and denoted by $\br(k)$.

   $\br(k)$ is a torsion group, and 
  the order of $P$ in $\br(k)$ is called the {\it period} of   $P$.

  \end{defn}
 
\begin{lem} \label{P.min.lem.2}  Let $\langle P\rangle \subset \br(k)$ denote the subgroup generated by $P$.  Then $Q\in \langle P\rangle$ iff there is a rational map $P\map Q$. 
  If these hold then $\ind(Q)$ divides $\ind(P)$.
  \qed
\end{lem}

\begin{cor} \label{BR.BC.COR} 
   Let $P$ be  a Severi--Brauer variety. Then the kernel of the base change map  $\br(k)\to \br\bigl(k(P)\bigr)$ is the subgroup generated by $P$. \qed
 \end{cor}

The stable birational equivalence version of Theorem~\ref{stab.bir.thm} is proved in \cite{MR70624}. 
 
 \begin{lem} \label{stab.bir.lem} 
 Let $P,  Q$ be  Severi--Brauer varieties   over a  field $k$.
 Then  $\langle P\rangle=\langle Q\rangle$ 
 iff   $P\times \p^{\dim Q}$ is birational to $Q\times\p^{\dim P} $.
     \qed
 \end{lem}

 \begin{lem}[Primary decomposition]\label{primary.lem}
   Let $P,  Q$ be  minimal Severi--Brauer varieties   over a  field $k$.
   Assume that  $(\dim P+1, \dim Q+1)=1$.
   Then $P\otimes Q$ is the minimal Severi--Brauer variety
   in the  Brauer equivalence class of $[P]\cdot [Q]$.\qed
   \end{lem}

\subsection*{Proof of Theorem~\ref{bir.of.prod.thm} $\Rightarrow$ 
  Theorem~\ref{stab.bir.thm}}{\ }

  \begin{say}\label{stab.bir.thm.pf}
    By Lemma~\ref{stab.bir.lem}, if  $P$ is birational to $Q$ then they  generate the same subgroup of $\br(k)$.
    The new result is the converse.

    Assume first that $P$ is not minimal.   By Lemma~\ref{P.min.lem}
    $$
    P\bire P^{\rm min}\times \p^r,\qtq{and}
    Q\bire Q^{\rm min}\times \p^r,
    $$
    where $r=\dim P-\dim P^{\rm min}\geq \dim P^{\rm min}+1$.
    Thus
    $$
    P^{\rm min}\times \p^r \bire Q^{\rm min}\times \p^r
    $$
    by Lemma~\ref{stab.bir.lem}, and we are done. (Here we did not need the assumption about  2 distinct prime factors.)

    We are left with the case when $P,Q$ are minimal. 
  By assumption we can write $\dim P+1=(m+1)(n+1)$ where $(m+1, n+1)=1$, and  $m, n$ are
  both $\geq 1$. Choose $a, c>0$ such that
  $$a(m+1)+c(n+1)\equiv 1  \mod  \dim P+1.
  $$
  Let $P_1$ (resp.\ $P_2$) be the minimal Severi--Brauer variety that is Brauer equivalent to $P^{a(m+1)}$ (resp.\ $P^{c(n+1)}$), and similarly for 
  $Q_1, Q_2$.  Note that
  $P\bire  P_1\times P_2$ and $Q\bire  Q_1\times Q_2$ by Lemma~\ref{primary.lem}. 

  Then $P_1, Q_1$ have index $n+1$, dimension $n$,  and generate the same subgroup of $\br(k)$. Similarly, $P_2, Q_2$ have index $m+1$, dimension $m$, and generate the same subgroup of $\br(k)$.
  Therefore, by Lemma~\ref{stab.bir.lem}, 
  $$
  P_1\times \p^{n}\bire  Q_1\times \p^{n},\qtq{and}
  P_2\times \p^{m}\bire  Q_2\times \p^{m}.
  $$
  By  (\ref{bir.of.prod.thm}.1)
  $$
  P\cong P_1\otimes P_2\bire  P_1\times P_2\times \p^{nm}.
   $$
  Since $mn\geq m,n$, we can first exchange  $P_1$ to $Q_1$, and then
  $P_2$ to $Q_2$ in the product. Thus we get that
  $$
  P\bire  P_1\times P_2\times \p^{nm}\bire  Q_1\times Q_2\times \p^{nm}
  \bire Q. \qed
  $$
  \end{say}

\subsection*{Roquette's method}{\ }

The geometric version is based on the following.

 \begin{lem}\label{span.lem} Let $L_0, \dots, L_n\subset \p^N$ be linear spaces. Assume that they span $\p^N$ and 
   $\sum (1+\dim L_i)=1+N$.  Let $p\in \p^N$ be a general point.

   Then there is a unique linear subvariety $M\subset \p^N$ of dimension $n$ that contains $p$ and meets every $L_i$.
 \end{lem}

 Proof. We use induction on $n$. If $n=0$, then $M=\{p\}$.
 For $n>0$, project from $L_n$ to get
 $L'_1, \dots, L'_{n-1}\subset \p^{N-r}$, where $r=1+\dim L_n$.
 Induction gives $M'\subset \p^{N-r}$ of dimension $n{-}1$.
 Then $M$ is the span of $M'$ and $p$. \qed

 \smallskip

 The proof shows that `$p\in \p^N$  general' means: not contained in the span of any 
 $n$ of the $L_i$.

 \medskip

 In order to state a version for Severi--Brauer varieties, we need  the notion of Weil restriction, introduced in Weil's 1959 lectures;  see  \cite[Sec.I.3]{MR670072}. It is called the {\it trace} in
 \cite{MR166215, MR2363485}.

 \begin{defn}[Weil restriction] \label{weil.rest} A general definition is given in \cite[Sec.7.6]{blr}. We use it only for geometrically normal varieties and separable field extension, for which the following simpler definition works.

   Let $K/k$ be a finite, separable field extension, and $X$ a
   geometrically normal $K$-variety. Its {\it Weil restriction}
   $\Re_{K/k}(X)$ is the unique, normal $k$-variety with the following property.

   \smallskip
  (\ref{weil.rest}.1)  Let $k'/k$ be a field extension and $K':=K\otimes_kk'$.
   Then $K'$-points of $X$ naturally correspond to
   $k'$-points of  $\Re_{K/k}(X)$.
   \smallskip

   One should think of  $\Re_{K/k}(X)$ as the product of the
   Galois conjugates of $X$ over $k$, with the natural Galois action. In particular, $\Re_{K/k}(X)$ becomes isomorphic to the product of $|K:k|$ copies of $X_{\bar k}$ over $\bar k$.

   The Weil restriction of a Severi--Brauer variety is not a Severi--Brauer variety, but
   $$
   \Re_{K/k}(\p^n_K)\bire \p^{r}_k, \qtq{where} r=n\cdot |K:k|.
   \eqno{(\ref{weil.rest}.2)}
   $$

 \end{defn}

 The key observation is the following.

 \begin{cor} \label{span.lem.sb.cor} Let $P$ be a Severi--Brauer variety over a field $k$. Let $K/k$ be a field extension of degree $n$, and $L\subset P_K$ a
   twisted linear subvariety. Assume that the conjugates $\{L^{\sigma}:\sigma\in G\}$ span $L$, and are linearly independent.

   Then there is a twisted linear subbundle
   $R_L(P)\subset \Re_{K/k}(L)\times P$
   such that the coordinate projection $R_L(P)\to P$ is birational.
   \end{cor}

 Proof.
 Let $k'/k$ be a field extension, $K':=K\otimes_kk'$, and 
 $p_0\in P_{K'}$  a general $K'$-point with $K'/k'$-conjugates
 $p_0,\dots, p_m$.

 By definition, a $k'$-point of $\Re_{K/k}(L)$ is the same as a
 $K'$-point $p_0\in P_{K'}$, which can be identified with
 the set of its $K'/k'$-conjugates
 $p_0,\dots, p_n$.

 Let 
 $R_L(P)\subset \Re_{K/k}(L)\times P$ be the
 twisted linear subbundle whose fiber over
 $p_0$ is the span of $p_0,\dots, p_n$.

 Birationality of the coordinate projection $R_L(P)\to P$ can be checked over $\bar k$, where Lemma~\ref{span.lem} applies. \qed

 \medskip

 Note that  $L$ and   the fibers of $R_L(P)\to \Re_{K/k}(L)$  have lower dimension than $P$. 
 The method of \cite{MR166215} uses  information about
 them to
get results about $P$.

\subsection*{Proof of Theorem~\ref{bir.of.prod.thm}}{\ }

We can get more  if we have a good  birational description of $\Re_{K/k}(L)$. The answer seems complicated in general, but simple in the case that we need.
The only new result is the statement (\ref{W.r.lem}.1), all the steps can be found in
\cite{MR166215, MR2363485}.

 \begin{prop}  \label{W.r.lem} Let $P$ be a Severi--Brauer variety of dimension $n$ over $k$, and
   $K/k$ a separable field extension of degree $m+1<n+1$.
Assume that $(m+1, n+1)=1$. 
   Then
   $$
   \Re_{K/k}(P_K)\bire P\times \p^{mn}.
   \eqno{(\ref{W.r.lem}.1)}
   $$
 \end{prop}

 Proof. Let $\grass(m, P)$ be the Grassmannian parametrizing $m$-dimensional linear subspaces  $L\subset P$, and  $\univ(m, P)\to \grass(m, P)$
 the  universal family parametrizing $m$-dimensional pointed linear subspaces  $p\in L\subset P$. Sending $(p, L)\to \{p\}$ shows that
 $\univ(m, P)$ is a $\grass(m{-}1, \p^{n-1})$-bundle over $P$. Thus
 $$
 \univ(m, P)\bire  P\times \p^{m(n-m)}.
 \eqno{(\ref{W.r.lem}.2)}
 $$
 Every   fiber $F$ of $\univ(m, P)\to \grass(m, P)$ is a  
 Severi--Brauer variety of dimension $m$. The canonical class of $F$ is
 $\o_F(-m{-}1)$, and the restriction of the canonical class of $P$ is
 $\o_F(-n{-}1)$. Thus if  $(m+1, n+1)=1$ then  $\o_F(1)$ is a line bundle on $F$, hence
 $$
 \univ(m, P)\bire \grass(m, P)\times \p^m.
 \eqno{(\ref{W.r.lem}.3)}
 $$
 Therefore
 $$
 \grass(m, P)\times \p^m\bire  P\times \p^{m(n-m)}.
 \eqno{(\ref{W.r.lem}.4)}
 $$
 After these preliminaries, let
$k'/k$ be a field extension, $K':=K\otimes_kk'$, and 
 $p_0\in P_{K'}$  a general $K'$-point with $K'/k'$-conjugates
 $p_0,\dots, p_m$. They span a linear subspace $L(p_0)$ of dimension $m$, which is defined over $k'$. We can thus first choose $L(p_0)$ in $\grass(m, P)(k')$, and then the point $p_0$ in $L(p_0)(K')$, which is  the corresponding   fiber of $\univ(m, P)\to \grass(m, P)$.
 Thus $\Re_{K/k}(P_K)$ is birational to the relative Weil restriction
 of $\univ(m, P)\to \grass(m, P)$. Using (\ref{W.r.lem}.3) and
  (\ref{weil.rest}.2) we get that
 $$
 \Re_{K/k}(P_K)\bire \grass(m, P)\times \p^{m(m+1)}\bire  P \times \p^{mn}.\qed
 $$

{\it Remark \ref{W.r.lem}.5.}  If $(m+1, n+1)\neq 1$, then there is no general formula for 
  $ \Re_{K/k}(P_K)$ that involves only $P$ and $m$.

  Indeed, if $n=m$ then  $ \Re_{K/k}(P_K)$ has a $k$-point iff $P$ has a $K$-point. If $k$ is a number field, the latter holds for some choices of $K$, but not for others.

 \begin{say}[Proof of Theorem~\ref{bir.of.prod.thm}]\label{prod.bir.thm.pf}

   Since $P_1^{\rm min}\otimes P_2^{\rm min}$ is  Brauer equivalent to  $P_1\otimes P_2$, we may assume that $P_1, P_2$ are both minimal.
We may also assume that $\dim P_2<\dim P_1$.

   Let $j:P_1\times P_2\into P_1\otimes P_2$ be the Segre embedding.
   Let $K/k$ be a splitting field of $P_2$ of degree $\dim P_2+1$, and
   $p\in (P_2)_K$ a point.

   Apply Corollary~\ref{span.lem.sb.cor}
   to   $L:=j\bigl(P_1\times\{p\}\bigr)$.
   We get a
   twisted linear subbundle
   $$
   R_L(P_1\otimes P_2)\subset \Re_{K/k}(L)\times (P_1\otimes P_2)
   \eqno{(\ref{prod.bir.thm.pf}.1)}
   $$
   such that the coordinate projection $R_L(P_1\otimes P_2)\to P_1\otimes P_2$ is birational.

   Proposition~\ref{W.r.lem} shows that
   $$
   \Re_{K/k}(L)\bire P_1\times \p^r, \qtq{where}
   r=\dim P_1\cdot \dim P_2.
   \eqno{(\ref{prod.bir.thm.pf}.2)}
   $$
   Using Corollary~\ref{BR.BC.COR},
   this implies that pulling $ P_1\otimes P_2$ back to $\Re_{K/k}(L)$ trivializes  $P_1$, thus $ P_1\otimes P_2$ becomes Brauer equivalent to $P_2$ over $\Re_{K/k}(L)$.
   Since  $(\dim P_1+1, \dim P_2+1)=1$, the index of $P_2$ over $\Re_{K/k}(L)$ is still $\dim P_2+1$. Thus in fact
   $$
   R_L(P_1\otimes P_2)\bire \Re_{K/k}(L)\times  P_2.
   $$
   Combining this with Corollary~\ref{span.lem.sb.cor}  and
   (\ref{prod.bir.thm.pf}.2) we get that 
   $$
   P_1\otimes P_2\bire R_L(P_1\otimes P_2)\bire \Re_{K/k}(L)\times  P_2
   \bire P_1\times  P_2 \times \p^r,
   $$
   where
   $r=\dim P_1\cdot \dim P_2$. \qed
   \end{say}

We have the following partial result without the relatively prime condition.

 \begin{lem}  \label{prod.to.power.lem}
   Let $P, Q$ be  Severi--Brauer varieties  over a field $k$. Then there are rational maps
   $P\otimes Q\map  P^{\dim Q+1}$.
 \end{lem}

 Proof. Geometrically, consider the Segre embedding
 $j:\p^n\times \p^m\into \p^N$. Pick a spanning set of points $p_0,\dots, p_n\in \p^n$ and set $L_i:= j\bigl(\{p_i\}\times \p^m\bigr)$. Pick  sections
 $$
 s_i\in H^0\bigl(L_i, \o_{\p^N}(r)|_{L_i}\bigr).
 $$
 Let $W_r\subset H^0\bigl(\p^N, \o_{\p^N}(r)\bigr)$ be the subspace of those sections whose restriction to $L_i$ is a multiple of $s_i$.
 This defines a rational map  $\p^N\map \p^M$ which
 induces the $r$-fold Veronese embedding on each copy of
 $\p^n\times \{q\}$ for $q\in \p^m$.

 Now let $P, Q$ be a Severi--Brauer varieties, $p_i$ a conjugate set of points and choose $r=\dim Q+1$. Although $ \o_{P\otimes Q}(\dim Q+1)$ is a twisted line bundle, its restriction to $L_i$ is the anticanonical line bundle, hence we can pick  $s_0$
 over the field $k(p_0)$ and let $s_i$ be its conjugates. Then the above construction gives a rational map
 $P\otimes Q\to  P^{\dim Q+1}$. \qed

 \begin{ack} I thank  J.-L.~Colliot-Th\'el\`ene, D.~Krashen  and 
T.~Szamuely  for 
comments, corrections  and references.
 Partial  financial support    was provided  by the Simons Foundation   (grant number SFI-MPS-MOV-00006719-02).
\end{ack}
 

\begin{thebibliography}{BLR90}

\bibitem[Alb39]{MR595}
A.~Adrian Albert, \emph{Structure of {A}lgebras}, American Mathematical Society
  Colloquium Publications, vol. Vol. 24, American Mathematical Society, New
  York, 1939. \MR{595}

\bibitem[Ami55]{MR70624}
S.~A. Amitsur, \emph{Generic splitting fields of central simple algebras}, Ann.
  of Math. (2) \textbf{62} (1955), 8--43. \MR{70624}

\bibitem[BLR90]{blr}
Siegfried Bosch, Werner L{\"u}tkebohmert, and Michel Raynaud, \emph{N\'eron
  models}, Ergebnisse der Mathematik und ihrer Grenzgebiete (3), vol.~21,
  Springer-Verlag, Berlin, 1990. \MR{1045822 (91i:14034)}

\bibitem[GS06]{MR2266528}
Philippe Gille and Tam{\'a}s {Sz}amuely, \emph{Central simple algebras and
  {G}alois cohomology}, Cambridge Studies in Advanced Mathematics, vol. 101,
  Cambridge University Press, Cambridge, 2006. \MR{2266528}

\bibitem[Her68]{MR227205}
I.~N. Herstein, \emph{Noncommutative rings}, The Carus Mathematical Monographs,
  vol. No. 15, Mathematical Association of America, ; distributed by John Wiley
  \& Sons, New York, 1968. \MR{227205}

\bibitem[Hog09]{MR2439423}
Amit Hogadi, \emph{Products of {B}rauer-{S}everi surfaces}, Proc. Amer. Math.
  Soc. \textbf{137} (2009), no.~1, 45--50. \MR{2439423}

\bibitem[Kol05]{k-conic}
J{\'a}nos Koll{\'a}r, \emph{Conics in the {G}rothendieck ring}, Adv. Math.
  \textbf{198} (2005), no.~1, 27--35. \MR{2183248 (2006k:14064)}

\bibitem[Kol25]{k-sb}
\bysame, \emph{{Severi-Brauer varieties; a geometric treatment}},
  \url{https://arxiv.org/abs/1606.04368}, 2025.

\bibitem[Kra03]{MR2029172}
Daniel Krashen, \emph{Severi-{B}rauer varieties of semidirect product
  algebras}, Doc. Math. \textbf{8} (2003), 527--546. \MR{2029172}

\bibitem[Kra08]{MR2363485}
\bysame, \emph{Birational maps between generalized {S}everi-{B}rauer
  varieties}, J. Pure Appl. Algebra \textbf{212} (2008), no.~4, 689--703.
  \MR{2363485}

\bibitem[Roq64]{MR166215}
Peter Roquette, \emph{Isomorphisms of generic splitting fields of simple
  algebras}, J. Reine Angew. Math. \textbf{214/215} (1964), 207--226.
  \MR{166215}

\bibitem[Ser62]{MR150130}
Jean-Pierre Serre, \emph{Corps locaux}, Publications de l'Institut de
  Math\'ematique de l'Universit\'e{} de Nancago, vol. VIII, Hermann, Paris,
  1962, Actualit\'es Scientifiques et Industrielles, No. 1296. \MR{150130}

\bibitem[Ser79]{MR554237}
\bysame, \emph{Local fields}, Graduate Texts in Mathematics, vol.~67,
  Springer-Verlag, New York-Berlin, 1979, Translated from the French by Marvin
  Jay Greenberg. \MR{554237}

\bibitem[Tre91]{tregub}
S.~L. Tregub, \emph{Birational equivalence of {B}rauer-{S}everi manifolds},
  Uspekhi Mat. Nauk \textbf{46} (1991), no.~6(282), 217--218. \MR{1164209}

\bibitem[Wei82]{MR670072}
Andr\'e Weil, \emph{Adeles and algebraic groups}, Progress in Mathematics,
  vol.~23, Birkh\"auser, Boston, MA, 1982, With appendices by M. Demazure and
  Takashi Ono. \MR{670072}

\end{thebibliography}

 \def\cprime{$'$} \def\cprime{$'$} \def\cprime{$'$} \def\cprime{$'$}
  \def\cprime{$'$} \def\dbar{\leavevmode\hbox to 0pt{\hskip.2ex
  \accent"16\hss}d} \def\cprime{$'$} \def\cprime{$'$}
  \def\polhk#1{\setbox0=\hbox{#1}{\ooalign{\hidewidth
  \lower1.5ex\hbox{`}\hidewidth\crcr\unhbox0}}} \def\cprime{$'$}
  \def\cprime{$'$} \def\cprime{$'$} \def\cprime{$'$}
  \def\polhk#1{\setbox0=\hbox{#1}{\ooalign{\hidewidth
  \lower1.5ex\hbox{`}\hidewidth\crcr\unhbox0}}} \def\cdprime{$''$}
  \def\cprime{$'$} \def\cprime{$'$} \def\cprime{$'$} \def\cprime{$'$}
\providecommand{\bysame}{\leavevmode\hbox to3em{\hrulefill}\thinspace}
\providecommand{\MR}{\relax\ifhmode\unskip\space\fi MR }
\providecommand{\MRhref}[2]{%
  \href{http://www.ams.org/mathscinet-getitem?mr=#1}{#2}
}
\providecommand{\href}[2]{#2}

 \bigskip

\noindent  Princeton University, Princeton NJ 08544-1000

{\begin{verbatim} kollar@math.princeton.edu\end{verbatim}}

\end{document}